\documentclass[11pt]{article}
\usepackage{amsfonts}
\usepackage{latexsym,amssymb,amsmath,graphics,cite}
\usepackage{tikz}
\usepackage{fancyhdr}
\usepackage{lipsum}
\usepackage{lineno,slashbox}

\begin{document}
\newcommand{\qed}{\hphantom{.}\hfill $\Box$\medbreak}
\newcommand{\Proof}{\noindent{\bf Proof \ }}
\newtheorem{theorem}{Theorem}[section]
\newtheorem{proposition}[theorem]{Proposition}
\newtheorem{lemma}[theorem]{Lemma}
\newtheorem{corollary}[theorem]{Corollary}
\newtheorem{remark}[theorem]{Remark}
\newtheorem{example}[theorem]{Example}
\newtheorem{definition}[theorem]{Definition}
\newtheorem{construction}[theorem]{Construction}
\renewcommand{\theequation}{\arabic{section}.\arabic{equation}}


\begin{center}
{\Large\bf Frame difference families and resolvable balanced incomplete block designs\footnote{Supported by NSFC under Grant $11471032$, and Fundamental Research Funds for the Central Universities under Grant $2016$JBM$071$, $2016$JBZ$012$ (T. Feng), NSFC under Grant 11771227, and Zhejiang Provincial Natural Science Foundation of China under Grant LY17A010008 (X. Wang).}}

\vskip12pt

Simone Costa$^a$, Tao Feng$^b$, Xiaomiao Wang$^c$\\[2ex]
{\footnotesize $^a$Dipartimento DICATAM, Universit\`a degli Studi di Brescia, Via Valotti 9, I-25123 Brescia, Italy}\\
{\footnotesize $^b$Department of Mathematics, Beijing Jiaotong University, Beijing 100044, P. R. China}\\
{\footnotesize $^c$Department of Mathematics, Ningbo University, Ningbo 315211, P. R. China}\\
{\footnotesize simone.costa@unibs.it, tfeng@bjtu.edu.cn, wangxiaomiao@nbu.edu.cn}
\vskip12pt

\end{center}

\vskip12pt

\noindent {\bf Abstract:} Frame difference families, which can be obtained via a careful use of cyclotomic conditions attached to strong difference families, play an important role in direct constructions for resolvable balanced incomplete block designs. We establish asymptotic existences for several classes of frame difference families. As corollaries new infinite families of 1-rotational $(pq+1,p+1,1)$-RBIBDs over $\mathbb{F}_{p}^+ \times \mathbb{F}_{q}^+$ are derived, and the existence of $(125q+1,6,1)$-RBIBDs is discussed. We construct $(v,8,1)$-RBIBDs for $v\in\{624,1576,2976,5720,5776,10200,14176,24480\}$, whose existence were previously in doubt. As applications, we establish asymptotic existences for an infinite family of optimal constant composition codes and an infinite family of strictly optimal frequency hopping sequences.




\section{Introduction}
Throughout this paper, sets and multisets will be denoted by curly braces $\{\ \}$ and square brackets $[\ ]$, respectively. When we emphasize a set or a multiset is fixed with an ordering, we regard it as a sequence, and denote it by $(\ )$. Every union will be understood as {\em multiset union} with multiplicities of elements preserved. $A{\cup}A{\cup}\cdots{\cup}A$ ($h$ times) will be denoted by $\underline{h}A$. If $A$ and $B$ are multisets defined on a multiplicative group, then $A\cdot B$ denotes the multiset $[ab:a\in A,b\in B]$.

For a positive integer $v$, we abbreviate $\{0,1,\dots,v-1\}$ by $\mathbb{Z}_v$ or $I_v$, with the former indicating that a cyclic group of this order is acting.

A $(v,k,\lambda)$-BIBD $(${\em balanced incomplete block design}$)$ is a pair $(V,\cal{A})$ where $V$ is a set of $v$ {\em points} and $\cal A$ is a collection of $k$-subsets of $X$ $($called {\em blocks}$)$ such that every
$2$-subset of $X$ is contained in exactly $\lambda$ blocks of $\cal A$. A $(v,k,\lambda)$-BIBD $(V,\cal{A})$ is said to be {\em resolvable}, or briefly a $(v,k,\lambda)$-RBIBD, if there exists a partition $\cal R$ of $\cal A$ $($called a {\em resolution}$)$ into {\em parallel classes}, each of which is a partition of $V$.

A powerful idea to obtain resolvable designs is given by the use of a special class of relative difference families: frame difference families. This concept was put forward by M. Buratti in 1999 \cite{b99}.

Let $(G,+)$ be an abelian group of order $g$ with a subgroup $N$ of order $n$. A $(G,N,k,\lambda)$ {\em relative difference family} $(DF)$, or $(g,n,k,\lambda)$-DF over $G$ relative to $N$, is a family $\mathfrak{B}=[B_1,B_2,\dots,B_r]$ of $k$-subsets of $G$ such that the list
$$\Delta \mathfrak{B}:=\bigcup_{i=1}^r[x-y:x,y\in B_i, x\not=y]=\underline{\lambda}(G\setminus N),$$
i.e., every element of $G\setminus N$ appears exactly $\lambda$ times in the multiset $\Delta \mathfrak{B}$ while it has no element of $N$. The members of $\mathfrak{B}$ are called {\em base blocks} and the number $r$ equals to $\lambda(g-n)/(k(k-1))$. A $(G,\{0\},k,\lambda)$-DF is said to be a {\em difference set} if it contains only one base block, written simply as $(G,k,\lambda)$-DS or $(g,k,\lambda)$-DS over $G$. The {\em complement} of a $(g,k,\lambda)$-DS over $G$ with the base block $B$ is a $(g,g-k,g-2k+\lambda)$-DS over $G$ with the base block $G\setminus B$.

Let $\mathfrak{F}$ be a $(g,n,k,\lambda)$-DF over $G$ relative to $N$. $\mathfrak{F}$ is a {\em frame difference family} $(FDF)$ if it can be partitioned into $\lambda n/(k-1)$ subfamilies $\mathfrak{F}_1,\mathfrak{F}_2,\ldots,\mathfrak{F}_{\lambda n/(k-1)}$ such that each $\mathfrak{F}_i$ has size of $(g-n)/(nk)$, and the union of base blocks in each $\mathfrak{F}_i$ is a system of representatives for the nontrivial cosets of $N$ in $G$. When $\lambda n=k-1$, a $(g,n,k,\lambda)$-FDF is said to be {\em elementary}.

The following proposition reveals the relation between frame difference families and resolvable designs, which can be seen as a corollary by combining the results of Theorem 1.1 in \cite{b99} and Theorem 5.11 in \cite{gm}. We outline the proof for completeness.

\begin{proposition}\label{FDFtoRBIBD}
If there exist a $(G,N,k,\lambda)$-FDF and a $(|N|+1,k,\lambda)$-RBIBD, then there exists a $(|G|+1,k,\lambda)$-RBIBD.
\end{proposition}

\Proof Let $\mathfrak{F}$ be a $(G,N,k,\lambda)$-FDF, which can be partitioned into $\lambda |N|/(k-1)$ subfamilies $\mathfrak{F}_1,\mathfrak{F}_2,\ldots,\mathfrak{F}_{\lambda |N|/(k-1)}$ such that $\bigcup_{F\in \mathfrak{F}_i, h\in N}(F+h)=G\setminus N$ for each $1\leq i\leq \lambda|N|/(k-1)$. Set $$\mathfrak{P}_i=\{F+h:F\in \mathfrak{F}_i, h\in N\}$$
for $1\leq i\leq \lambda|N|/(k-1)$. Let $S$ be a complete system of representatives for the cosets of $N$ in $G$. For each $s\in S$, construct a $(|N|+1,k,\lambda)$-RBIBD on $(N+s)\cup \{\infty\}$, where $\infty\not\in G$. It has $\lambda |N|/(k-1)$ parallel classes, written as $\mathfrak{Q}_{s,i}$, $1\leq i\leq \lambda|N|/(k-1)$. It is readily checked that $(\mathfrak{P}_i+s)\cup \mathfrak{Q}_{s,i}$, $1\leq i\leq \lambda|N|/(k-1)$ and $s\in S$, constitute all parallel classes of a $(|G|+1,k,\lambda)$-RBIBD, which is defined on $G\cup \{\infty\}$. \qed

An {\em automorphism group} of a $(v,k,\lambda)$-RBIBD $(V,{\cal A})$ with $\cal R$ as its resolution is a group of permutations on $V$ leaving $\cal A$ and $\cal R$ invariant, respectively. A $(v,k,\lambda)$-RBIBD is said to be {\em $1$-rotational} over a group $G$ of order $v-1$ if it admits $G$ as an automorphism group fixing one point and acting sharply transitively on the others.

By revisiting the proof of Proposition \ref{FDFtoRBIBD}, one can have the following proposition.

\begin{proposition}\label{FDFtoRBIBD 1-rotational}
Suppose there exists a $(G,N,k,\lambda)$-FDF. If there is a $1$-rotational $(|N|+1,k,\lambda)$-RBIBD over $N$, then there is a $1$-rotational $(|G|+1,k,\lambda)$-RBIBD over $G$.
\end{proposition}

The target of this paper is to construct frame difference families and (1-rotational) resolvable designs via strong difference families. Let $\mathfrak{S}=[F_1,F_2,\dots,F_s]$ with $F_i=(f_{i,0},f_{i,1},\dots,f_{i,k-1})$ for $1\leq i\leq s$, be a family of $s$ multisets of size $k$ defined on a group $(G,+)$ of order $g$. We say that $\mathfrak{S}$ is a $(G,k,\mu)$ {\em strong difference family}, or a $(g,k,\mu)$-SDF over $G$, if the list
$$\Delta \mathfrak{S}:=\bigcup_{i=1}^s [f_{i,a}-f_{i,b}: 0\leq a,b\leq k-1, a\not=b]=\underline{\mu} G,$$ i.e., every element of $G$ (0 included) appears exactly $\mu$ times in the multiset $\Delta \mathfrak{S}$. The members of $\mathfrak{S}$ are called {\em base blocks} and the number $s$ equals to $\mu g/(k(k-1))$. Note that $\mu$ is necessarily even since the element $0\in G$ is expressed in even ways as differences in any multiset.

\begin{proposition}\label{prop:nece SDF}
A $(G,k,\mu)$-SDF exists only if $\mu$ is even and $\mu |G|\equiv 0 \pmod{k(k-1)}$.
\end{proposition}

The concept of strong difference families was introduced in \cite{b99} and revisited in \cite{bg,m}. It is useful in the constructions of relative difference families and BIBDs (cf. \cite{cfw}), and perfect cycle decompositions (cf. \cite{bcw}). Many direct constructions for RBIBDs in the literature are obtained by the use of certain suitable SDFs explicitly or implicitly (cf. \cite{b97,b99,bf,bz}).

Although many authors have worked on the existence of $(v,8,1)$-RBIBDs, there are still $66$ open cases for small values of $v$ (see Table 4 in \cite{ga} or Table 7.41 in \cite{ag}). In Section 2 we shall show that, with a careful application of cyclotomic conditions attached to a strong difference family, we can establish the existence of three new $(v,8,1)$-RBIBDs for $v\in\{624,1576,2976\}$. Then via known recursive constructions for RBIBDs, we obtain another five new $(v,8,1)$-RBIBDs for $v\in\{5720,5776,10200,14176,24480\}$.

M. Buratti, J. Yan and C. Wang \cite{byw} proved that any $(k-1,k,kt)$-SDF can lead to a $((k-1)p,k-1,k,1)$-FDF for any sufficiently large prime $p$ and $p\equiv kt+1 \pmod{2kt}$. We shall generalize their result in Section 3 (see Theorem \ref{thm:FDF-2}).

In Sections 4 and 5, we shall prove that, if the initial SDF has some particular patterns, then the lower bound on $q$ can be reduced greatly. Theorem \ref{thm:DF-4} generalizes Construction A in \cite{bf}, and Theorems \ref{thm:DF-3} and \ref{thm:Singer DS-2} generalize Construction B in \cite{bf}. As corollaries of Theorems \ref{thm:DF-3}-\ref{thm:Singer DS-2} and \ref{thm:DF-4}, Theorems \ref{thm:BIBD-3}-\ref{thm:BIBD-4} give new 1-rotational $(pq+1,p+1,1)$-RBIBDs. Theorem \ref{6,1RBIBD2} presents a new infinite family of $(v,6,1)$-RBIBDs.

As applications, in Section 6, we derive new optimal constant composition codes and new strictly optimal frequency hopping sequences.

\section{Basic lemma and new $(v,8,1)$-RBIBDs}

Let $q$ be a prime power. As usual we denote by $\mathbb{F}_q$ the finite field of order $q$, by $\mathbb{F}_q^+$ its additive group, by $\mathbb{F}^*_q$ its multiplicative group, by $\mathbb{F}_{q}^{\Box}$ the set of nonzero squares, and by $\mathbb{F}_{q}^{\not\Box}$ nonsquares in $\mathbb{F}_{q}$.

If $q\equiv 1 \pmod{e}$, then $C_0^{e,q}$ will denote the group of nonzero $e$th powers of $\mathbb{F}_q$ and once a primitive element $\omega$ of $\mathbb{F}_q$ has been fixed, we set $C_i^{e,q}=\omega^i\cdot C_0^{e,q}$ for $i=0,1,\ldots,e-1$. We refer to the cosets $C_0^{e,q},C_1^{e,q},\ldots,C_{e-1}^{e,q}$ of $C_0^{e,q}$ in $\mathbb{F}^*_q$ as the {\em cyclotomic classes of index $e$}. Let $A$ be a multisubset of $\mathbb{F}_q^*$. If each cyclotomic coset $C_l^{e,q}$ for $l\in I_e$ contains exactly $\lambda$ elements of $A$, then $A$ is said to be a {\em $\lambda$-transversal} for these cosets. If $A$ is a $1$-transversal, $A$ is often referred to as a {\em representative system for the cosets} of $C_0^{e,q}$ in $\mathbb{F}_q^*$. The following lemma allows us to obtain frame different families by using strong difference families.

\begin{lemma}\label{lem:FDF-2}
Let $q\equiv 1 \pmod {e}$ be a prime power and $d|e$. Let $S$ be a representative system for the cosets of $C_0^{e,q}$ in $C_0^{d,q}$. Let $d(q-1)\equiv 0\pmod{ek}$ and $t=d(q-1)/ek$. Suppose that there exists a $(G,k,kt\lambda)$-$SDF$ $\mathfrak{S}=[F_1,F_2,\ldots, F_n]$, where $\lambda|G|\equiv 0\pmod{k-1}$ and $F_i=(f_{i,0},f_{i,1},\dots,f_{i,k-1})$, $1\leq i\leq n$. If there exists a partition  $\mathcal{P}$ of base blocks of $\mathfrak{S}$ into $\lambda|G|/(k-1)$ multisets, each of size $t$, such that one can choose appropriate multiset $[\Phi_1,\Phi_2,\dots,\Phi_n]$ of ordered $k$-subsets of $\mathbb{F}_q^*$ with $\Phi_i=(\phi_{i,0},\phi_{i,1},\dots,\phi_{i,k-1})$, $1\leq i \leq n$, satisfying that
\begin{itemize}
\item[$(1)$] $\bigcup_{i=1}^n [\phi_{i,a}-\phi_{i,b}: f_{i,a}-f_{i,b}=h, (a,b)\in I_k\times I_k, a\not=b]=C_0^{e,q}\cdot D_h$ for each $h\in G$, where $D_h$ is a $\lambda$-transversal for the cosets of $C_0^{d,q}$ in $\mathbb{F}_q^*$,
\item[$(2)$] $\bigcup_{i:F_i\in P} [\phi_{i,a}: a\in I_k]=C_0^{e,q}\cdot E_P$ for each $P \in\mathcal{P}$, where $E_P$ is a representative system for the cosets of $C_0^{d,q}$ in $\mathbb{F}_q^*$,
\end{itemize}
then
$$ \mathfrak{F}=[B_i\cdot \{(1,s)\}: 1\leq i\leq n, s\in S]$$
is a $(G\times \mathbb{F}_q^+,G\times \{0\},k,\lambda)$-FDF, where $B_i=\{(f_{i,0},\phi_{i,0}),(f_{i,1},\phi_{i,1}),\dots,(f_{i,k-1},\phi_{i,k-1})\}$.
\end{lemma}

\Proof Since $n=\lambda kt|G|/(k(k-1))$ and $t=d(q-1)/ek$, we have
$$|\mathfrak{F}|= n\cdot \frac{e}{d}=\frac{\lambda|G|(q-1)}{k(k-1)},$$
which coincides with the number of base blocks of a $(G\times \mathbb{F}_q^+,G\times \{0\},k,\lambda)$-FDF. According to $\mathcal{P}$, we can partition $\mathfrak{F}$ into $\lambda|G|/(k-1)$ subfamilies
$$\mathfrak{F}_P=[B_i\cdot \{(1,s)\}: F_i\in P, s\in S],$$
where $P\in \mathcal{P}$. Each subfamily contains $|P|\times|S|=te/d=(q-1)/k$ base blocks. Because of Condition (2),
$$\bigcup_{i: F_i\in P}\bigcup_{s\in S} B_i\cdot\{(1,s)\}$$
forms a representative system for the nontrivial cosets of $G\times \{0\}$ in $G\times \mathbb{F}_q^+$. Finally it is readily checked that
\begin{eqnarray*}
\Delta \mathfrak{F}=\hspace{-3mm}&\hspace{-3mm}& \bigcup_{s\in S} \bigcup^n_{i=1}(\Delta B_i\cdot\{(1,s)\})=\bigcup_{s\in S} \bigcup^n_{i=1} [(f_{i,a}-f_{i,b},(\phi_{i,a}-\phi_{i,b})\cdot s): (a,b)\in I_k\times I_k, a\not=b]\\
=\hspace{-3mm}&\hspace{-3mm}&\bigcup_{s\in S}[ \{h\}\times (C_0^{e,q}\cdot D_h\cdot \{s\}): h\in G ]={\underline\lambda} (G\times \mathbb{F}_q^*).
\end{eqnarray*}
Therefore $\mathfrak{F}$ is a $(G\times \mathbb{F}_q^+,G\times \{0\},k,\lambda)$-FDF. \qed

\subsection{A 1-rotational $(624,8,1)$-RBIBD}
In this subsection we shall apply Lemma \ref{lem:FDF-2} with $e=q-1$ to present a new RBIBD. Note that when $e=q-1$, $C_0^{e,q}=\{1\}$ and $S=C_0^{d,q}$.

\begin{lemma}\label{FDF1}
There exists an elementary $(\mathbb{Z}_{7}\times \mathbb{F}_{89}^+,\mathbb{Z}_{7}\times \{0\},8,1)$-FDF.
\end{lemma}

\Proof Take the $(\mathbb{Z}_{7},8,8)$-SDF containing the unique base block $[0,0,1,1,2,2,4,4]$ as the first components of base blocks of the required FDF. Let
$$B=\{(0,1),(0,20),(1,14),(1,58),(2,18),(2,61),(4,26),(4,73)\}.$$
Then applying Lemma \ref{lem:FDF-2} with $G=\mathbb{Z}_{7}$, $q=89$, $e=88$, $d=8$, $k=8$ and $\lambda=1$ which yield $|\mathcal{P}|=1$ and $t=1$, we have
$$\mathfrak{F}=[B\cdot (1,s):s\in C_0^{8,89}]$$
 forms an elementary $(\mathbb{Z}_{7}\times \mathbb{F}_{89}^+,\mathbb{Z}_{7}\times \{0\},8,1)$-FDF. It is readily checked that each $D_h$, $h\in\mathbb{Z}_{7}$, is a representative system for the cosets of $C_0^{8,89}$ in $\mathbb{F}_{89}^*$ (for example $D_0=\{19,42,43,44,45,46,47,70\}$ and $D_1=\{3,4,13,38,47,49,57,83\}$). The unique $E_P=\{1$, $20,14,58,18,61,26,73\}$ is also a representative system for the cosets of $C_0^{8,89}$ in $\mathbb{F}_{89}^*$. \qed

\begin{theorem}\label{624}
There exists a 1-rotational $(624,8,1)$-RBIBD over $\mathbb{Z}_{623}$.
\end{theorem}

\Proof By Lemma \ref{FDF1} there exists a $(\mathbb{Z}_{7}\times \mathbb{F}_{89}^+,\mathbb{Z}_{7}\times \{0\},8,1)$-FDF. Apply Proposition \ref{FDFtoRBIBD 1-rotational} with a trivial 1-rotational $(8,8,1)$-RBIBD to obtain a 1-rotational $(624,8,1)$-RBIBD over $\mathbb{Z}_{7}\times \mathbb{F}_{89}^+$ that is isomorphic to $\mathbb{Z}_{623}$. \qed

\subsection{A $(v,8,1)$-RBIBD  for $v\in\{1576,2976\}$}
In this subsection we shall apply Lemma \ref{lem:FDF-2} with $e=\frac{q-1}{4}$ to present two new RBIBDs.
\begin{lemma}\label{lem:SDF(63,8,8)}
There exists a $(\mathbb{Z}_{p},8,8)$-SDF for $p\in \{63,119\}$.
\end{lemma}

\Proof For $p=63$, take \begin{center}
\begin{tabular}{lll}
$F_1=[20,20,-20,-20,29,29,-29,-29]$, \\
$F_2=F_3=F_4=F_5=[0,1,3,7,19,34,42,53]$,\\
$F_6=F_7=F_8=F_9=[0,1,4,6,26,36,43,51]$.
\end{tabular}
\end{center}
\noindent Then the multiset $[F_i: 1\leq i\leq 9]$ forms a $(\mathbb{Z}_{63},8,8)$-SDF.

For $p=119$, take \begin{center}
\begin{tabular}{lll}
$F_1=[20,20,-20,-20,29,29,-29,-29],$\\
$F_2=F_3=F_4=F_5=[0,1,42,28,101,97,94,114],$\\
$F_6=F_7=F_8=F_9=[0,1,12,23,41,85,104,106],$\\
$F_{10}=F_{11}=F_{12}=F_{13}=[0,2,5,17,37,47,68,76],$\\
$F_{14}=F_{15}=F_{16}=F_{17}=[0,4,10,38,54,62,86,93].$
\end{tabular}
\end{center}
\noindent Then the multiset $[F_i: 1\leq i\leq 17]$ forms a $(\mathbb{Z}_{119},8,8)$-SDF.
\qed

\begin{lemma}\label{lem:63 25}
There exists a $(\mathbb{Z}_{p}\times \mathbb{F}_{25}^+,\mathbb{Z}_{p}\times \{0\},8,1)$-FDF for $p\in \{63,119\}$.
\end{lemma}
\Proof Take the $(\mathbb{Z}_{p},8,8)$-SDF from Lemma \ref{lem:SDF(63,8,8)} as the first components of base blocks of the required $(\mathbb{Z}_{p}\times\mathbb{F}_{25}^+,\mathbb{Z}_{p}\times\{0\},8,1)$-FDF. Take $x^2-x+2$ to be a primitive polynomial of degree $2$ over $\mathbb{F}_{5}$ and $\omega$ to be a primitive root in $\mathbb{F}_{25}$. Let $\xi=\omega^6$.
For $p=63$, let
\begin{center}
\begin{tabular}{l}
$B_1= \{(20,1),(20,-1),(-20,\xi),(-20,-\xi),(29, \omega),(29,-\omega),(-29,\omega\xi),(-29,-\omega\xi)\},$\\
$B_2= \{(0,1),(1,\omega^{17}),(3,\omega^{12}),(7,\omega^{5}),(19,\omega^{23}),(34,\omega^{11}),(42,\omega^{18}),(53,\omega^{6})\},$\\
$B_6= \{(0,1),(1,\omega^{11}),(4,\omega^{5}),(6,\omega^{12}),(26,\omega^{23}),(36,\omega^{18}),(43,\omega^{17}),(51,\omega^{6})\},$\\
$B_3= B_2\cdot\{(1,-1)\},$\ \ \ \ \ \ \ \ \ \ \ $B_4=B_2\cdot\{(1,\xi)\}$,\ \ \ \ \ \ \ \ \ \ \ $B_5= B_2\cdot\{(1,-\xi)\}$,\\
$B_7=B_6\cdot\{(1,-1)\}$,\ \ \ \ \ \ \ \ \ \ \
$B_8= B_6\cdot\{(1,\xi)\},$\ \ \ \ \ \ \ \ \ \ \ $B_9=B_6\cdot\{(1,-\xi)\}$.\\
\end{tabular}
\end{center}
For $p=119$, let
\begin{center}
\begin{tabular}{l}
$B_1=\{(20,1),(20,-1),(-20,\xi),(-20,-\xi),(29,\omega),(29,-\omega),(-29,\omega\xi),(-29,-\omega\xi)\},$\\
$B_2=\{(0,\omega),(1,\omega^7),(42,1),(28,\omega^6),(101,\omega^{18}),(97,\omega^{13}),(94,\omega^{12}),(114,\omega^{19})\},$\\
$B_6=\{(0,1),(1,\omega^{12}),(12,\omega),(23,\omega^{18}),(41,\omega^{13}),(85,\omega^7),(104,\omega^6),(106,\omega^{19})\},$\\
$B_{10}=\{(0,\omega),(2,\omega^7),(5,\omega^{12}),(17,\omega^6),(37,\omega^{19}),(47,\omega^{18}),(68,1),(76,\omega^{13})\},$\\
$B_{14}=\{(0,\omega),(4,\omega^4),(10,\omega^{10}),(38,\omega^7),(54,\omega^{22}),(62,\omega^{19}),(86,\omega^{16}),(93,\omega^{13})\}.$\\
$B_3= B_2\cdot\{(1,-1)\},$\ \ \ \ \ \ \ \ \ \ \ $B_4=B_2\cdot\{(1,\xi)\}$,\ \ \ \ \ \ \ \ \ \ \ $B_5= B_2\cdot\{(1,-\xi)\}$,\\
$B_7=B_6\cdot\{(1,-1)\}$,\ \ \ \ \ \ \ \ \ \ \
$B_8= B_6\cdot\{(1,\xi)\},$\ \ \ \ \ \ \ \ \ \ \ $B_9=B_6\cdot\{(1,-\xi)\}$,\\
$B_{11}= B_{10}\cdot\{(1,-1)\},$\ \ \ \ \ \ \ \ \ $B_{12}=B_{10}\cdot\{(1,\xi)\}$,\ \ \ \ \ \ \ \ \ $B_{13}= B_{10}\cdot\{(1,-\xi)\}$,\\
$B_{15}=B_{14}\cdot\{(1,-1)\}$,\ \ \ \ \  \ \ \ \
$B_{16}=B_{14}\cdot\{(1,\xi)\},$\ \ \ \ \ \ \ \ \ $B_{17}=B_{14}\cdot\{(1,-\xi)\}$.\\
\end{tabular}
\end{center}
Let $S$ be a representative system for the cosets of $C_0^{6,25}=\{1,-1,\xi,-\xi\}$ in $C_0^{2,25}$. Then, applying Lemma \ref{lem:FDF-2} with $G=\mathbb{Z}_{p}$, $q=25$, $e=6$, $d=2$, $k=8$ and $\lambda=1$ which yield $|\mathcal{P}|=p/7$ and $t=1$, we have a $(\mathbb{Z}_{p}\times \mathbb{F}_{25}^+,\mathbb{Z}_{p}\times \{0\},8,1)$-FDF for $p\in \{63,119\}$. Note that for any $i$, $1 \leq i \leq p/7$, $\bigcup_{s\in S} B_i\cdot\{(1,s)\}$
forms a representative system for the nontrivial cosets of $\mathbb{Z}_{p}\times \{0\}$ in $\mathbb{Z}_{p}\times \mathbb{F}_{25}^+$. \qed

\begin{theorem}\label{thm:RBIBD(1576)}
There exists a $(1576,8,1)$-RBIBD and a $(2976,8,1)$-RBIBD.
\end{theorem}
\Proof By Lemma \ref{lem:63 25}, there exists a $(\mathbb{Z}_{p}\times \mathbb{F}_{25}^+,Z_{p}\times \{0\},8,1)$-FDF for $p\in \{63,119\}$. Apply Proposition \ref{FDFtoRBIBD} with a $(p+1,8,1)$-RBIBD, which exists by Table 7.41 in \cite{ag}, to obtain a $(1576,8,1)$-RBIBD and a $(2976,8,1)$-RBIBD. \qed

\subsection{Five new RBIBDs via recursive constructions}

A {\em transversal design} is a triple $(V,\mathcal{G},\mathcal{B})$, where $V$ is a set of $km$ points, $\mathcal{G}$ is a partition of $V$ into $k$ {\em groups}, each of size $m$, and $\mathcal{B}$ is a set of $k$-subsets $($called {\em blocks}$)$ of $V$ satisfying every pair of $V$ is contained either in exactly one group or in exactly one block, but not both. Such a design is denoted by a TD$(k,m)$.

It is well known that the existence of a TD$(k,m)$ is equivalent to the existence of $k-2$ mutually orthogonal Latin squares of order $m$. A TD$(q+1,q)$ exists for any prime power $q$ (cf. Theorem 6.44 in \cite{StinsonBook}), and a TD$(10,48)$ exists by Theorem 2.1 in \cite{ac}.

\begin{lemma}\label{aGre}{\rm (Lemma 4.9 in \cite{ga})}
Suppose there exist a TD$(10,m)$ and a $(56m+8,8,1)$-RBIBD. For any given $0\leq n\leq m$, if there exists a $(56n+8,8,1)$-RBIBD, then there exists a $(56(9m+n)+8,8,1)$-RBIBD.
\end{lemma}

\begin{lemma}\label{aGreBis}{\rm (Lemma 4.34 in \cite{ga})}
Suppose there exist a TD$(9,8n)$, a $(56n+8,8,1)$-RBIBD and a $(56m+8,8,1)$-RBIBD. Then there exists a $(56(8mn+n)+8,8,1)$-RBIBD.
\end{lemma}

\begin{theorem}\label{5776}
There exists a $(v,8,1)$-RBIBD for $v\in\{5720,5776,10200,24480\}$.
\end{theorem}

\Proof
Apply Lemma \ref{aGreBis} with $m=2$ and $n=6$ to obtain a $(5720,8,1)$-RBIBD, where the needed $(v,8,1)$-RBIBDs for $v\in\{120,344\}$ are from Table 7.41 in \cite{ag}.

For $v\in\{5776,10200,24480\}$, apply Lemma \ref{aGre} with $(m,n)\in\{(11,4),(19,11),(48,5)\}$ to obtain a $(v,8,1)$-RBIBD, where the needed $(624,8,1)$-RBIBD is from Theorem \ref{624} and the needed $(u,8,1)$-RBIBDs for $u\in\{232,288,1072,2696\}$ are from Table 7.41 in \cite{ag}. \qed

Let $(H,+)$ be an abelian group of order $h$. An $(H,k,\lambda)$ {\em difference matrix} $($briefly, $(H,k,\lambda)$-DM$)$ is a $k\times h\lambda$ matrix $D=(d_{ij})$ with entries from $H$ so that for each $1\leq i<j\leq k$ the multiset
$$\{ d_{il}-d_{jl}: 1\leq l\leq h\lambda \}$$
contains every element of $H$ exactly $\lambda$ times. An $(H,k,1)$-DM is {\em homogeneous} if its each row is a permutation of elements of $H$.

The property of a difference matrix is preserved even if one add any element of $H$ to all entries in any row or column of the difference matrix. Then, w.l.o.g., all entries in the first row in a difference matrix are zero. Such a difference matrix is said to be {\em normalized}. Any normalized difference matrix can yield a homogeneous difference matrix by deleting its first row. Thus the existence of a homogeneous $(H,k-1,1)$-DM is equivalent to that of an $(H,k,1)$-DM. The multiplication table for the finite field $\mathbb{F}_q$ is an $(\mathbb{F}_q^+,q,1)$-DM \cite{d}.

The following construction is a variation of standard recursive construction for difference families (see for example Theorem 6.1 in \cite{byw}).

\begin{construction}\label{recursive}
Suppose there exists a $(G,N,k,1)$-FDF. If there exists a homogeneous $(H,k,1)$-DM, then there exists a $(G\times H,N\times H,k,1)$-FDF.
\end{construction}

\begin{theorem}\label{thm:RBIBD(14176)}
There exists a $(14176,8,1)$-RBIBD.
\end{theorem}

\Proof Take a $(\mathbb{Z}_{63}\times \mathbb{F}_{25}^+,Z_{63}\times \{0\},8,1)$-FDF from Lemma \ref{lem:63 25}. Then apply Construction \ref{recursive} with a homogeneous $(\mathbb{F}_{9}^+,8,1)$-DM to obtain a $(\mathbb{Z}_{63}\times \mathbb{F}_{25}^+\times \mathbb{F}_{9}^+,Z_{63}\times \{0\}\times \mathbb{F}_{9}^+,8,1)$-FDF. Finally apply Proposition \ref{FDFtoRBIBD} with a $(568,8,1)$-RBIBD, which exists by Table 7.41 in \cite{ag}, to obtain a $(14176,8,1)$-RBIBD. \qed

\section{Asymptotic existence of FDFs}

Throughout this paper we always write
\begin{eqnarray*}\label{Q(e,m)}
Q(d,m)=\frac{1}{4}(U+\sqrt{U^2+4d^{m-1}m})^2, \mbox{ where } U=\sum_{h=1}^m {m \choose h}(d-1)^h(h-1)
\end{eqnarray*}
for given positive integers $d$ and $m$. The following theorem characterizes existences of elements satisfying certain cyclotomic conditions in a finite field.

\begin{theorem}\label{thm:cyclot bound} {\rm \cite{bp,cj}}
Let $q\equiv 1 \pmod{d}$ be a prime power, let $B=\{b_0,b_1,\dots, b_{m-1}\}$ be an arbitrary m-subset of $\mathbb{F}_q$ and let $(\beta_0,\beta_1,\dots,\beta_{m-1})$ be an arbitrary element of $\mathbb{Z}_d^m$. Set $X=\{x\in \mathbb{F}_q: x-b_i\in C^{d,q}_{\beta_i} \mbox{ for } i=0,1,\dots,m-1\}$. Then $X$ is not empty for any prime power $q\equiv 1 \pmod{d}$ and $q>Q(d,m)$.
\end{theorem}

Abel and Buratti \cite{ab0} announced Theorem \ref{thm:cyclot bound} without proving it. The case of $m=3$ in Theorem \ref{thm:cyclot bound} was first shown by Buratti \cite{b02}. Then a proof similar to that of $m=3$ allows Chang and Ji \cite{cj}, and Buratti and Pasotti \cite{bp} to generalize this result to any $m$. Theorem \ref{thm:cyclot bound} is derived from Weil's Theorem (see \cite{ln}, Theorem 5.41) on multiplicative character sums and plays an essential role in the asymptotic existence problem for difference families (cf. \cite{cwz}).

The idea of the following lemma is from Theorem 4.1 in \cite{byw}.

\begin{lemma}\label{thm:FDF-1}
If there exists a $(G,k,kt\lambda)$-SDF with $\lambda|G|\equiv 0\pmod{k-1}$, then there exists a $(G\times \mathbb{F}_q^+,G\times \{0\},k,\lambda)$-FDF
\begin{itemize}
\item for any even $\lambda$ and any prime power $q\equiv 1 \pmod{kt}$ with $q>Q(kt,k)$;
\item for any odd $\lambda$ and any prime power $q\equiv kt+1 \pmod{2kt}$ with $q>Q(kt,k)$.
\end{itemize}
\end{lemma}

\Proof By assumption one can take a $(G,k,kt\lambda)$-SDF $\mathfrak{S}=[F_{1},F_{2},\dots,F_{n}]$, where $n=t\lambda|G|/(k-1)$ and $F_i=[f_{i,0},f_{i,1},\ldots,f_{i,k-1}]$, $1\leq i\leq n$. To apply Lemma \ref{lem:FDF-2} with $e=q-1$ and $d=kt$, we need to give a partition $\mathcal{P}$ of base blocks of $\mathfrak{S}$ into $\lambda|G|/(k-1)$ multisets, each of size $t$, such that one can choose appropriate multiset $[\Phi_1,\Phi_2,\dots,\Phi_n]$ of ordered $k$-subsets of $\mathbb{F}_q^*$ with $\Phi_i=(\phi_{i,0},\phi_{i,1},\dots,\phi_{i,k-1})$, $1\leq i \leq n$, satisfying that
\begin{itemize}
\item[$(1)$] $\bigcup_{i=1}^n [\phi_{i,a}-\phi_{i,b}: f_{i,a}-f_{i,b}=h, (a,b)\in I_k\times I_k, a\not=b]=D_h$ for each $h\in G$,
\item[$(2)$] $\bigcup_{i: F_i\in P} [\phi_{i,a}: a\in I_k]= E_P$ for each $P \in\mathcal{P}$,
\end{itemize}
where $D_h$ is a $\lambda$-transversal for the cosets of $C_0^{d,q}$ in $\mathbb{F}_q^*$ and $E_P$ is a representative system for the cosets of $C_0^{d,q}$ in $\mathbb{F}_q^*$.

Thus if one can choose an appropriate mapping $\pi$ acting on symbolic expressions satisfying that
\begin{itemize}
\item[$(1')$] $\bigcup_{i=1}^n [\pi(\phi_{i,a}-\phi_{i,b}): f_{i,a}-f_{i,b}=h, (a,b)\in I_k\times I_k, a\not=b]=\underline{\lambda}I_d$ for each $h\in G$,
\item[$(2')$] $\bigcup_{i:F_i\in P} [\pi(\phi_{i,a}): a\in I_k]=I_d$ for each $P \in\mathcal{P}$.
\end{itemize}
and can choose appropriate elements of $\Phi_i$, $1\leq i\leq n$, such that these elements are consistent with the mapping $\pi$, i.e., $\pi$ can be seen as a function from $\mathbb{F}_q^*$ to $\mathbb{Z}_d$ satisfying $\pi(x)=\theta$ if $x\in C^{d,q}_\theta$, then one can apply Lemma \ref{lem:FDF-2} to obtain a $(G\times \mathbb{F}_q^+,G\times \{0\},k,\lambda)$-FDF.

This procedure can always be done. First we need to choose an appropriate mapping $\pi$ satisfying Conditions $(1')$ and $(2')$. Condition $(2')$ can be satisfied easily. The key is how to meet Condition $(1')$. Let $G_2$ denote the subgroup of $\{h\in G:2h=0\}$. When $\lambda$ is even, we can specify $\pi$ to satisfy
$$\left\{
\begin{array}{ll}
\bigcup_{i=1}^n [\pi(\phi_{i,a}-\phi_{i,b}): f_{i,a}-f_{i,b}=h, (a,b)\in I_k\times I_k, a\not=b]=\underline{\lambda}I_d, & h\in G\setminus G_2,\\\\
\bigcup_{i=1}^n [\pi(\phi_{i,a}-\phi_{i,b}): f_{i,a}-f_{i,b}=h, (a,b)\in I_k\times I_k, a<b]=\underline{\frac{\lambda}{2}} I_d, & h\in G_2.\\
\end{array}
\right.
$$
When $\lambda$ is odd, we can specify $\pi$ to satisfy
$$\left\{
\begin{array}{ll}
\bigcup_{i=1}^n [\pi(\phi_{i,a}-\phi_{i,b}): f_{i,a}-f_{i,b}=h, (a,b)\in I_k\times I_k, a\not=b]=\underline{\lambda}I_d, & h\in G\setminus G_2,\\\\
\bigcup_{i=1}^n [\pi(\phi_{i,a}-\phi_{i,b}): f_{i,a}-f_{i,b}=h, (a,b)\in I_k\times I_k, a<b]={\underline \lambda} I_{\frac{d}{2}}, & h\in G_2.\\
\end{array}
\right.
$$
Note that when $\lambda$ is odd, by Proposition \ref{prop:nece SDF}, the existence of the given $(G,k,kt\lambda)$-SDF implies $d=kt$ is even. When $q\equiv d+1 \pmod{2d}$, $-1\in C_{\frac{d}{2}}^{d,q}$.

Once $\pi$ is fixed, one can apply Theorem \ref{thm:cyclot bound} and Lemma \ref{lem:FDF-2} to obtain the required $(G\times \mathbb{F}_q^+,G\times \{0\},k,\lambda)$-FDF. \qed

\begin{theorem}\label{thm:FDF-2}
If there exists a $(G,k,kt\lambda)$-SDF with $\lambda|G|\equiv 0\pmod{k-1}$, then there exists a $(G\times \mathbb{F}_q^+,G\times \{0\},k,\lambda)$-FDF
\begin{itemize}
\item for any even $\lambda$ and any prime power $q\equiv 1 \pmod{kt}$ with $q>Q(kt,k)$;
\item for any odd $\lambda$ and any prime power $q\equiv krt+1 \pmod{2krt}$ with $q>Q(krt,k)$, where $r$ is any positive integer.
\end{itemize}
\end{theorem}

\Proof For even $\lambda$, the conclusion is straightforward by Lemma \ref{thm:FDF-1}. For odd $\lambda$, if a $(G,k,kt\lambda)$-SDF exists with $n$ base blocks, then there exists a $(G,k,krt\lambda)$-SDF with $rn$ base blocks for any positive integer $r$, which means by Lemma \ref{thm:FDF-1} that there exists a $(G\times \mathbb{F}_q^+,G\times \{0\},k,\lambda)$-FDF for any prime power $q\equiv krt+1 \pmod{2krt}$ with $q>Q(krt,k)$. \qed

\section{FDFs from SDFs with particular patterns}

The application of Theorem \ref{thm:FDF-2} results in a huge lower bound on $q$. To reduce the lower bound, we shall request the initial SDF has some special patterns. If an SDF only contains one base block, then it is referred to as a {\em difference multiset} (cf. \cite{b99}) or a {\em regular difference cover} (cf. \cite{abms}).

\begin{lemma}\label{lem:SDF-Paley}{\rm \cite{b99}}
\begin{itemize}
\item[$(1)$] Let $p\equiv 3\pmod{4}$ be a prime power. Then ${\underline 2}(\{0\}\cup\mathbb{F}^{\Box}_p)$ is an $(\mathbb{F}_p^+,p+1,p+1)$-SDF $($called Paley difference multiset of the second type$)$.

\item[$(2)$] Let $p$ be an odd prime power. Set $X_1={\underline 2}(\{0\}\cup\mathbb{F}^{\Box}_p)$ and $X_2={\underline 2}(\{0\}\cup\mathbb{F}^{\not \Box}_p)$. Then $[X_1,X_2]$ is an $(\mathbb{F}_p^+,p+1,2p+2)$-SDF $($called Paley strong difference family of the third type$)$.

\item[$(3)$] Given twin prime powers $p>2$ and $p+2$, the set
$(\mathbb{F}^{\Box}_p\times \mathbb{F}^{\Box}_{p+2})\cup (\mathbb{F}^{\not\Box}_p\times \mathbb{F}^{\not\Box}_{p+2})\cup(\mathbb{F}_p\times\{0\})$ is a $(p(p+2),\frac{p(p+2)-1}{2},\frac{p(p+2)-3}{4})$-DS over $\mathbb{F}_p^+\times \mathbb{F}_{p+2}^+$. Let $D$ be its complement. Then ${\underline 2}D$ is a $(p(p+2),p(p+2)+1,p(p+2)+1)$ difference multiset $($called twin prime power difference multiset$)$.

\item[$(4)$] Given any prime power $p$ and any integer $m\geq3$, there is a $(\frac{p^m-1}{p-1},\frac{p^{m-1}-1}{p-1},\frac{p^{m-2}-1}{p-1})$ difference set over $\mathbb{Z}_{\frac{p^m-1}{p-1}}$. Let $D$ be its complement. Then ${\underline p}D$ is a $(\frac{p^m-1}{p-1},p^m,p^m(p-1))$ difference multiset $($called Singer difference multiset$)$.
\end{itemize}
\end{lemma}

By Theorem \ref{thm:FDF-2} we can easily obtain an infinite family of FDFs from each of the SDFs in Lemma \ref{lem:SDF-Paley}. For example, by the second type Paley $(\mathbb{F}_p^+,p+1,p+1)$-SDF, we get

\begin{corollary}\label{cor:Paley DS-1}
Let $p\equiv 3\pmod{4}$ be a prime power. Then there exists an $(\mathbb{F}_{p}^+ \times \mathbb{F}_{q}^+,\mathbb{F}_{p}^+\times \{0\},p+1,1)$-FDF for any prime power $q\equiv p+2 \pmod{2(p+1)}$ and $q>Q(p+1,p+1)$.
\end{corollary}

The lower bound on $q$ in Corollary \ref{cor:Paley DS-1} is huge even if $p$ is small. For example, if $p=11$, then $Q(12,12)=7.94968\times 10^{27}$. Thus it would be meaningful to develop a new technique to reduce the bound.

\begin{lemma}\label{lem:DF-3}
Let $G$ be an additive group of odd order $l$. Suppose that there exists a $(G,l+1,l+1)$-SDF whose unique base block $(f_{0},f_1,\ldots,f_{l})=$
\begin{eqnarray*}\label{2st Paley}
(x_0,x_0,x_1,x_1,\ldots,x_{\frac{l-1}{2}},x_{\frac{l-1}{2}}),
\end{eqnarray*}
where $x_0,x_1,\ldots,x_{(l-1)/2}$ are distinct elements of $G$.
Let $q$ be a prime power satisfying $q\equiv 1\pmod{l+1}$ and let  $d=(l+1)/2$ .
Suppose that one can choose an appropriate multiset $(\phi_{0},\phi_1,\ldots,\phi_{l})=$
\begin{eqnarray*}\label{2st Paley-y}
(y_0,-y_0,y_1,-y_1,\ldots,y_{\frac{l-1}{2}},-y_{\frac{l-1}{2}})
\end{eqnarray*}
such that $\{y_0,y_1,\ldots,y_{(l-1)/2}\}\subseteq {\mathbb F}_{q}^*$ and for each $h\in G$,
\begin{eqnarray*}\label{2st Paley-dh}
[\phi_{a}-\phi_{b}:f_{a}-f_{b}=h,(a,b)\in I_{l+1}\times I_{l+1},a\neq b]=\{1,-1\}\cdot D_h,
\end{eqnarray*}
where $D_h$ is a representative system for the cosets of $C_0^{d,q}$ in $\mathbb{F}_{q}^*$. Let $S$ be a representative system for the cosets of $\{1,-1\}$ in $C_0^{d,q}$. Let $B=\{(f_0,\phi_0),(f_1,\phi_1),\ldots,(f_{l},\phi_{l})\}.$
Then
$$\mathfrak{F}=[B\cdot\{(1,s)\}:s\in S]$$
forms an elementary $(G\times \mathbb{F}_{q}^+,G\times \{0\},l+1,1)$-FDF.
\end{lemma}

\Proof Since $d$ is a divisor of $l+1$ and $q\equiv 1\pmod{l+1}$, $d$ is also a divisor of $q-1$. This makes $C_0^{d,q}$ meaningful. The assumption $q\equiv 1\pmod{l+1}$ ensures $-1\in C_0^{d,q}$. Since $q$ is odd, $y_i\neq -y_i$ for any $0\leq i\leq (l-1)/2$, $B\cdot\{(1,s)\}$ is a set of size $l+1$ for any $s\in S$. Then applying Lemma \ref{lem:FDF-2} with $e=(q-1)/2$, $d=(l+1)/2$, $k=l+1$ and $\lambda=1$ which yield $|\mathcal{P}|=1$ and $t=1$, we have a $(G\times \mathbb{F}_{q}^+,G\times \{0\},l+1,1)$-FDF. Note that $\{1,-1\}=C_0^{(q-1)/2,q}$ and $\{y_0,-y_0,y_1,-y_1,\ldots,y_{\frac{l-1}{2}},-y_{\frac{l-1}{2}}\}=\{1,-1\}\cdot \frac{1}{2}\cdot D_0$. \qed

\begin{lemma}\label{lem:DF-3-Dh}
Follow the notation in Lemma $\ref{lem:DF-3}$.
\begin{itemize}
\item[$(1)$] W.l.o.g., $D_0=2\cdot\{y_0,y_1,\ldots,y_{(l-1)/2}\}$.
\item[$(2)$]
For each $h\in G\setminus\{0\}$, let
$$T_h=[\phi_{a}-\phi_{b}:f_{a}-f_{b}=h,(a,b)\in I_{l+1}\times I_{l+1},a\neq b].$$
Then $T_h=\{1,-1\}\cdot D_h$ for some $D_h\subset \mathbb{F}_{q}$ and the size of $D_h$ is $(l+1)/2$. Furthermore, $D_h=D_{-h}$ and w.l.o.g., $D_h$ consists of elements of type $y_i\pm y_j$.
\item[$(3)$] Let $T$ be a representative system for the cosets of $\{1,-1\}$ in $G\setminus\{0\}$. Any element of type $y_i\pm y_j$ must be contained in a unique $D_h$ for some $h\in T$ $($note that the term ``element'' here is a symbolic expression; for example $y_1-y_2$ and $y_3+y_4$ are different element but they may have the same value$)$.
\end{itemize}\end{lemma}

\Proof The verification is straightforward. \qed

Begin with SDFs from Lemma \ref{lem:SDF-Paley}(1), (3) and (4), and then apply Lemma \ref{lem:DF-3}, \ref{lem:DF-3-Dh}, and Theorem \ref{thm:cyclot bound}. We have the following theorems. Note that we take $p=2$ when use Lemma \ref{lem:SDF-Paley}(4).

\begin{theorem}\label{thm:DF-3}
Let $p\equiv 3\pmod{4}$ be a prime power. Then there exists an elementary $(\mathbb{F}_{p}^+ \times \mathbb{F}_{q}^+,\mathbb{F}_{p}^+\times \{0\},p+1,1)$-FDF for any prime power $q\equiv 1 \pmod{p+1}$ and $q>Q((p+1)/2,p)$.
\end{theorem}

\begin{theorem}\label{thm:Twin DS-2}
Let $p$ and $p+2$ be twin prime powers satisfying $p>2$. Then there exists an elementary $(\mathbb{F}_{p}^+ \times \mathbb{F}_{p+2}^+\times\mathbb{F}_{q}^+,\mathbb{F}_{p}^+ \times \mathbb{F}_{p+2}^+\times \{0\},p(p+2)+1,1)$-FDF for any prime power $q\equiv 1 \pmod{p(p+2)+1}$ and $q>Q((p(p+2)+1)/2,p(p+2))$.
\end{theorem}

\begin{theorem}\label{thm:Singer DS-2}
Let $m\geq3$ be an integer. Then there exists an elementary $(\mathbb{Z}_{2^m-1} \times \mathbb{F}_{q}^+,\mathbb{Z}_{2^m-1} \times \{0\},2^m,1)$-FDF for any prime power $q\equiv 1 \pmod{2^{m}}$ and $q>Q(2^{m-1},2^{m}-1)$.
\end{theorem}

Compared with Corollary \ref{cor:Paley DS-1}, Theorem \ref{thm:DF-3} not only reduce the lower bound on $q$ but also relax the congruence condition on $q$. Take for example $p=11$. By Theorem \ref{thm:DF-3} the bound is $Q(6,11)=8.77844\times 10^{18}$, which is much smaller than $Q(12,12)$.

Theorems \ref{thm:DF-3} and \ref{thm:Singer DS-2} can be seen as a generalization of M. Buratti and N. Finizio's construction in \cite{bf} for $(\mathbb{Z}_{7} \times \mathbb{Z}_{q},\mathbb{Z}_{7}\times \{0\},8,1)$-FDFs, where $q$ is a prime.

\begin{lemma}\label{lem:DF-4}
Let $p\equiv 1\pmod{4}$ be a prime power and $d=p+1$. Let $q$ be a prime power satisfying $q\equiv 1\pmod{2p+2}$. Let $\omega$ be a generator of $\mathbb{F}_{p}$. Take the third type Paley $(\mathbb{F}_{p}^+,p+1,2p+2)$-SDF from Lemma $\ref{lem:SDF-Paley}(2)$ whose base blocks are:
\begin{center}
\begin{tabular}{l}
$(f_{10},f_{11},\ldots,f_{1p})=(0,0,\omega^2,\omega^2,\omega^4,\omega^4,\ldots,\omega^{p-1},\omega^{p-1}),$\\
$(f_{20},f_{21},\ldots,f_{2p})=(0,0,\omega,\omega,\omega^3,\omega^3,\ldots,\omega^{p-2},\omega^{p-2}).$
\end{tabular}
\end{center}
Suppose that one can choose appropriate multisets
\begin{center}
\begin{tabular}{l}
$(\phi_{10},\phi_{11},\ldots,\phi_{1p})=(y_0,-y_0,y_1,-y_1,y_2,-y_2,\ldots,y_{\frac{p-1}{2}},-y_{\frac{p-1}{2}}),$\\
$(\phi_{20},\phi_{21},\ldots,\phi_{2p})=(y_{\frac{p+1}{2}},-y_{\frac{p+1}{2}},y_{\frac{p+3}{2}},-y_{\frac{p+3}{2}},\ldots,y_{p},-y_{p}),$
\end{tabular}
\end{center}
such that $\{y_0,y_1,\ldots,y_{p}\}\subseteq {\mathbb F}_{q}^*$ and for each $h\in \mathbb{F}_{p}$,
\begin{eqnarray*}\label{2st Paley-dh}
\bigcup_{i=1}^2[\phi_{ia}-\phi_{ib}:f_{ia}-f_{ib}=h,(a,b)\in I_{p+1}\times I_{p+1},a\neq b]=\{1,-1\}\cdot D_h,
\end{eqnarray*}
where $D_h$ is a representative system for the cosets of $C_0^{d,q}$ in $\mathbb{F}_{q}^*$. Let $S$ be a representative system for the cosets of $\{1,-1\}$ in $C_0^{d,q}$. Let $B_i=\{(f_{i0},\phi_{i0}),(f_{i1},\phi_{i1}),\ldots,(f_{ip},\phi_{ip})\},$
where $i=1,2$. Then
$$\mathfrak{F}=\bigcup_{i=1}^2[B_i\cdot\{(1,s)\}:s\in S]$$
forms an elementary $(\mathbb{F}_{p}^+\times \mathbb{F}_{q}^+,\mathbb{F}_{p}^+\times \{0\},p+1,1)$-FDF.
\end{lemma}

\Proof Applying Lemma \ref{lem:FDF-2} with $G=\mathbb{F}_{p}^+$, $e=(q-1)/2$, $d=p+1$, $k=p+1$ and $\lambda=1$ which yield $|\mathcal{P}|=1$ and $t=2$, we obtain a  $(\mathbb{F}_{p}^+\times \mathbb{F}_{q}^+,\mathbb{F}_{p}^+\times \{0\},p+1,1)$-FDF. Note that $\{1,-1\}=C_0^{(q-1)/2,q}$ and $\{y_0,-y_0,y_1,-y_1,\ldots,y_p,-y_p\}=\{1,-1\}\cdot \frac{1}{2}\cdot D_0$. \qed

\begin{lemma}\label{lem:DF-4-Dh}
Follow the notation in Lemma $\ref{lem:DF-4}$.
\begin{itemize}
\item[$(1)$] W.l.o.g., $D_0=2\cdot\{y_0,y_1,y_2,\ldots,y_{p}\}$.
\item[$(2)$]
For each $h\in \mathbb{F}_{p}^*$, let
$$T_h=\bigcup_{i=1}^2[\phi_{ia}-\phi_{ib}:f_{ia}-f_{ib}=h,(a,b)\in I_{p+1}\times I_{p+1},a\neq b].$$
Then $T_h=\{1,-1\}\cdot D_h$ for some $D_h\subset \mathbb{F}_{q}$ and the size of $D_h$ is $p+1$. Furthermore, $D_h=D_{-h}$ and w.l.o.g., $D_h$ consists of elements of type $y_i\pm{y_j}$.
\item[$(3)$] Let $T$ be a representative system for the cosets of $\{1,-1\}$ in $\mathbb{F}_{p}^*$. Any element of type $y_i\pm{y_j}$ must be contained in a unique $D_h$ for some $h\in T$.
\end{itemize}\end{lemma}

\Proof The verification is straightforward. \qed

Combining the results of Lemmas \ref{lem:DF-4} and \ref{lem:DF-4-Dh}, and then applying Theorem \ref{thm:cyclot bound}, we have

\begin{theorem}\label{thm:DF-4}
Let $p\equiv 1\pmod{4}$ be a prime power. Then there exists an elementary $(\mathbb{F}_{p}^+ \times \mathbb{F}_{q}^+,\mathbb{F}_{p}^+\times \{0\},p+1,1)$-FDF for any prime power $q\equiv 1 \pmod{2p+2}$ and $q>Q(p+1,p)$.
\end{theorem}

Theorem \ref{thm:DF-4} can be seen as a generalization of M. Buratti and N. Finizio's construction in \cite{bf} for $(\mathbb{Z}_{5} \times \mathbb{Z}_{q},\mathbb{Z}_{5}\times \{0\},6,1)$-FDFs, where $q$ is a prime.

Start from the FDFs in Theorems \ref{thm:DF-3}, \ref{thm:Twin DS-2}, \ref{thm:Singer DS-2} and \ref{thm:DF-4}. Then apply Proposition \ref{FDFtoRBIBD 1-rotational} with a trivial 1-rotational $(k,k,1)$-RBIBD. We obtain the following theorems.

\begin{theorem}\label{thm:BIBD-3}
There exists a 1-rotational $(pq+1,p+1,1)$-RBIBD over $\mathbb{F}_{p}^+ \times \mathbb{F}_{q}^+$ for any prime powers $p$ and $q$ with $p\equiv 3 \pmod{4}$, $q\equiv 1\pmod{p+1}$ and $q>Q((p+1)/2,p)$.
\end{theorem}

\begin{theorem}\label{thm:Twin BIBD-2}
Let $p$ and $p+2$ be twin prime powers satisfying $p>2$. There exists a 1-rotational $(p(p+2)q+1,p(p+2)+1,1)$-RBIBD over $\mathbb{F}_{p}^+ \times\mathbb{F}_{p+2}^+ \times \mathbb{F}_{q}^+$ for any prime power $q\equiv 1 \pmod{p(p+2)+1}$ and $q>Q((p(p+2)+1)/2,p(p+2))$.
\end{theorem}

\begin{theorem}\label{thm:Singer BIBD-2}
There exists a 1-rotational $((2^m-1)q+1,2^m,1)$-RBIBD over $\mathbb{Z}_{2^m-1} \times \mathbb{F}_{q}^+$ for any integer $m\geq3$ and any prime power $q\equiv 1\pmod{2^m}$ and $q>Q(2^{m-1},2^m-1)$.
\end{theorem}

\begin{theorem}\label{thm:BIBD-4}
There exists a 1-rotational $(pq+1,p+1,1)$-RBIBD over $\mathbb{F}_{p}^+ \times \mathbb{F}_{q}^+$ for any prime powers $p$ and $q$ with $p\equiv 1 \pmod{4}$, $q\equiv 1\pmod{2p+2}$ and $q>Q(p+1,p)$.
\end{theorem}

\section{A family of RBIBDs with block size $6$}

\begin{lemma}\label{125SDF}
There exists a $(\mathbb{Z}_{125},6,6)$-SDF.
\end{lemma}

\Proof Take
\begin{center}
\begin{tabular}{ll}
$A_1=[0,0,19,19,71,71]$,&
$A_2=A_3=[0,10,28,51,78,97]$,\\
$A_4=A_5=[0,3,62,75,86,110]$, & $A_6=A_7=[0,5,12,58,70,112],$\\
$A_8=A_{9}=[0,7,27,44,70,96]$, & $A_{10}=A_{11}=[0,1,42,93,85,45],$\\
$A_{12}=A_{13}=[0,1,100,104,109,88]$, & $A_{14}=A_{15}=[0,1,90,81,21,32]$,\\
$A_{16}=A_{17}=[0,3,16,40,46,50]$, & $A_{18}=A_{19}=[0,2,7,29,35,68],$\\
$A_{20}=A_{21}=[0,2,8,57,102,116]$, & $A_{22}=A_{23}=[0,2,22,32,36,96],$\\
$A_{24}=A_{25}=[0,8,23,38,72,86].$ &
\end{tabular}
\end{center}
Then the multiset $[A_i: 1\leq i \leq 25]$ forms a $(\mathbb{Z}_{125}, 6,6)$-SDF. \qed

Applying Theorem \ref{thm:FDF-2} with a $(\mathbb{Z}_{125},6,6)$-SDF, we can obtain a $(\mathbb{Z}_{125}\times \mathbb{F}_q^+,\mathbb{Z}_{125}\times \{0\},6,1)$-FDF for any prime power $q\equiv 7 \pmod{12}$ and $q>Q(6,6)$. But the bound $Q(6,6)=3.4829\times 10^{10}$ is a little big. We shall reduce the bound by supplying a refined construction in this section.

\begin{lemma}\label{125FDF}
Let $q\equiv 7 \pmod{12}$ be a prime power. Take the $(\mathbb{Z}_{125},6,6)$-SDF given in Lemma $\ref{125SDF}$, whose base blocks are $A_1,A_2,\ldots,A_{25}$. Suppose one can choose appropriate multisets
$$C_1=(y_{11},-y_{11},y_{12},-y_{12},y_{13},-y_{13}),$$
$$C_{2i}=(y_{2i,1},-y_{2i,1},y_{2i,2},-y_{2i,2},y_{2i,3},y_{2i,4}),\hspace{0.5 cm}C_{2i+1}=-C_{2i},$$
where $1\leq i\leq 12$, such that each $C_j$, $1\leq j\leq 25$, is a representative system for the cosets of $C_0^{6,q}$ in $\mathbb{F}_q^*$. Furthermore, write $A_j=(a_{j1},a_{j2},a_{j3},a_{j4},a_{j5},a_{j6})$, $C_j=(c_{j1},c_{j2},c_{j3},c_{j4},c_{j5},c_{j6})$, and
$$B_j=\{(a_{j1},c_{j1}),(a_{j2},c_{j2}),(a_{j3},c_{j3}),(a_{j4},c_{j4}),(a_{j5},c_{j5}),(a_{j6},c_{j6})\},$$
where $1\leq j\leq 25$. Set $$\bigcup_{j=1}^{25} \Delta B_j=\bigcup_{l\in \mathbb{Z}_{125}} \{l\}\times \Delta_l.$$ If for any $l\in \mathbb{Z}_{125}$, $\Delta_l$  is a representative system for the cosets of $C_0^{6,q}$ in $\mathbb{F}_q^*$, then
$$\mathfrak{F}=[B_j\cdot \{(1,\alpha)\}: \alpha\in C_0^{6,q}, 1\leq j\leq 25]$$
forms a $(\mathbb{Z}_{125}\times \mathbb{F}_{q}^+,\mathbb{Z}_{125}\times \{0\},6,1)$-FDF.
\end{lemma}

\Proof Apply Lemma \ref{lem:FDF-2} with $G=\mathbb{Z}_{125}$, $e=q-1$, $d=6$, $k=6$ and $\lambda=1$ which yield $|\mathcal{P}|=25$ and $t=1$ to obtain the required  a $(\mathbb{Z}_{125}\times \mathbb{F}_{q}^+,\mathbb{Z}_{125}\times \{0\},6,1)$-FDF. \qed

\begin{lemma}\label{lem:125p}
Follow the notation in Lemma $\ref{125FDF}$.
\begin{itemize}
\item[$(1)$] $C_1=\{1,-1\}\cdot \{y_{11},y_{12},y_{13}\}.$
\item[$(2)$] For $1\leq i\leq 12$, $C_{2i}=(\{1,-1\}\cdot\{y_{2i,1},y_{2i,2}\})\cup\{y_{2i,3},y_{2i,4}\}$ and $C_{2i+1}=(\{1,-1\}\cdot \{y_{2i,1},y_{2i,2}\})\cup\{-y_{2i,3},-y_{2i,4}\}$.
\item[$(3)$] For each $l\in \mathbb{Z}_{125}$, $\Delta_l=\Delta_{-l}=\{1,-1\}\cdot D_l$ for some $D_l\subset \mathbb{F}_q$ and the size of $D_l$ is $3$.
\item[$(4)$] W.l.o.g., $D_0=2\cdot\{y_{11},y_{12},y_{13}\}$.
\item[$(5)$] For each $l\in \mathbb{Z}_{125}^*$, w.l.o.g., $D_l$ consists of elements having the following types
    \begin{center}
    \begin{tabular}{l}
    $(I)$ $y_{1,t_1}\pm y_{1,t_2}$ for some $1\leq t_1\neq t_2\leq 3$; \\
    $(II)$ $y_{2i,r}\pm y_{2i,s}$ for some $1\leq i\leq 12$, $r\in\{3,4\}$ and $s\in\{1,2\}$;\\
    $(III)$ $y_{2i,4}\pm y_{2i,3}$ for some $1\leq i\leq 12$;\\
    $(IV)$ $y_{2i,2}\pm y_{2i,1}$ for some  $1\leq i\leq 12$; \\
    $(V)$ $2y_{2i,s}$ for some $1\leq i\leq 12$ and $s\in\{1,2\}$.
    \end{tabular}
    \end{center}
\item[$(6)$] Any element of Types $(I)$, $(II)$, $(III)$ and $(V)$ is contained in a unique $D_l$ for some $1\leq l\leq 62$.
\item[$(7)$] Any element of Type $(IV)$ is contained in exactly two different $D_l$'s, say $D_{l_1}$ and $D_{l_2}$, for some $1\leq l_1\neq l_2\leq 62$.
\end{itemize}
\end{lemma}

\Proof It is readily checked that $(1)$-$(6)$ hold. The verification for $(7)$ is a little more complicated, which relies heavily on the given $(\mathbb{Z}_{125},6,6)$-SDF. For example, since $B_2=\{(0,y_{21}),(10,-y_{21}),(28,y_{22}),(51,-y_{22}),(78,y_{23}),(97,y_{24})\}$, we have $y_{22}-y_{21}\in D_{28}$ and $D_{41}$; $y_{22}+y_{21}\in D_{18}$ and $D_{51}$. By tedious calculation, one can check $(7)$. For convenience, we list each $D_l$ explicitly for $0\leq l\leq 62$ in Table \ref{tab3}. \qed

\begin{table}[t]{\tabcolsep 0.05in
{\small
\begin{tabular}{ll}
$D_0=[2y_{1,1},2y_{1,2},2y_{1,3}],$&
$D_1=[2y_{10,1},2y_{12,1},2y_{14,1}],$\\
$D_2=[2y_{18,1},2y_{20,1},2y_{22,1}],$&
$D_3=[2y_{4,1},y_{10,4}-y_{10,2},2y_{16,1}],$\\
$D_4=[2y_{12,2},y_{16,4}-y_{16,3},y_{22,3}+y_{22,2}],$ &
$D_5=[2y_{6,1},y_{12,3}+y_{12,2},y_{18,2}+y_{18,1}],$\\
$D_6=[y_{16,3}+y_{16,2},y_{18,3}+y_{18,2},y_{20,2}+y_{20,1}],$ &
$D_7=[y_{6,2}+y_{6,1},2y_{8,1},y_{18,2}-y_{18,1}],$\\
$D_8=[y_{10,3}+y_{10,2},y_{20,2}-y_{20,1},2y_{24,1}],$ &
$D_9=[y_{12,3}-y_{12,2},2y_{14,2},y_{20,4}-y_{20,1}],$\\
$D_{10}=[2y_{2,1},y_{16,4}+y_{16,2},2y_{22,2}],$ &
$D_{11}=[y_{4,3}+y_{4,2},y_{14,4}-y_{14,3},y_{20,4}+y_{20,1}],$\\
$D_{12}=[y_{6,2}-y_{6,1},y_{6,3}+y_{6,2},y_{12,4}-y_{12,2}],$ &
$D_{13}=[2y_{4,2},y_{6,4}-y_{6,1},y_{16,2}+y_{16,1}],$\\
$D_{14}=[y_{20,4}-y_{20,3},y_{22,3}-y_{22,2},y_{24,4}-y_{24,3}],$ &
$D_{15}=[y_{4,4}-y_{4,1},2y_{24,2},y_{24,2}+y_{24,1}],$\\
$D_{16}=[y_{12,3}-y_{12,1},y_{12,4}+y_{12,2},y_{16,2}-y_{16,1}],$ &
$D_{17}=[2y_{8,2},y_{12,3}+y_{12,1},y_{20,4}-y_{20,2}],$\\
$D_{18}=[y_{2,2}+y_{2,1},y_{4,4}+y_{4,1},y_{6,4}+y_{6,1}],$ &
$D_{19}=[y_{1,2}+y_{1,1},y_{1,2}-y_{1,1},y_{2,4}-y_{2,3}],$\\
$D_{20}=[y_{8,2}+y_{8,1},y_{14,3}+y_{14,1},y_{22,2}+y_{22,1}],$ &
$D_{21}=[y_{12,2}+y_{12,1},y_{12,4}-y_{12,3},y_{14,3}-y_{14,1}],$\\
$D_{22}=[y_{12,2}-y_{12,1},2y_{18,2},y_{22,2}-y_{22,1}],$ &
$D_{23}=[2y_{2,2},y_{20,3}-y_{20,1},y_{24,2}-y_{24,1}],$\\
$D_{24}=[y_{4,3}-y_{4,2},y_{4,4}-y_{4,3},2y_{16,2}],$ &
$D_{25}=[y_{6,4}-y_{6,2},y_{12,2}-y_{12,1},y_{20,3}+y_{20,1}],$\\
$D_{26}=[y_{8,3}+y_{8,2},y_{8,4}-y_{8,3},y_{12,2}+y_{12,1}],$ &
$D_{27}=[y_{2,3}+y_{2,2},y_{8,2}-y_{8,1},y_{18,2}-y_{18,1}],$\\
$D_{28}=[y_{2,2}-y_{2,1},y_{2,4}-y_{2,1},y_{18,3}-y_{18,2}],$ &
$D_{29}=[y_{8,4}-y_{8,1},y_{18,2}+y_{18,1},y_{22,4}-y_{22,1}],$\\
$D_{30}=[y_{16,3}-y_{16,2},y_{22,2}-y_{22,1},y_{24,2}-y_{24,1}],$ &
$D_{31}=[y_{14,4}+y_{14,1},y_{20,3}-y_{20,2},y_{22,4}+y_{22,1}],$\\
$D_{32}=[y_{10,2}+y_{10,1},y_{14,4}-y_{14,1},y_{22,2}+y_{22,1}],$ &
$D_{33}=[y_{10,2}-y_{10,1},y_{18,3}+y_{18,1},y_{18,4}-y_{18,3}],$\\
$D_{34}=[y_{16,4}-y_{16,2},y_{22,3}+y_{22,1},y_{24,3}+y_{24,2}],$ &
$D_{35}=[y_{4,4}+y_{4,2},y_{14,2}-y_{14,1},y_{18,3}-y_{18,1}],$\\
$D_{36}=[y_{8,4}+y_{8,1},y_{14,2}+y_{14,1},y_{22,3}-y_{22,1}],$ &
$D_{37}=[y_{8,2}-y_{8,1},y_{12,4}-y_{12,1},y_{16,2}-y_{16,1}],$\\
$D_{38}=[y_{2,4}+y_{2,1},y_{12,4}+y_{12,1},y_{24,2}+y_{24,1}],$ &
$D_{39}=[y_{4,3}-y_{4,1},y_{18,4}+y_{18,2},y_{24,4}-y_{24,1}],$\\
$D_{40}=[y_{10,3}-y_{10,1},y_{10,4}-y_{10,3},y_{16,2}+y_{16,1}],$ &
$D_{41}=[y_{2,2}-y_{2,1},y_{10,2}+y_{10,1},y_{10,3}+y_{10,1}],$\\
$D_{42}=[y_{4,3}+y_{4,1},y_{6,4}-y_{6,3},y_{10,2}-y_{10,1}],$ &
$D_{43}=[y_{8,3}-y_{8,2},y_{10,3}-y_{10,2},y_{16,3}+y_{16,1}],$\\
$D_{44}=[y_{8,2}+y_{8,1},y_{10,4}+y_{10,1},y_{14,2}+y_{14,1}],$ &
$D_{45}=[y_{10,4}-y_{10,1},y_{14,2}-y_{14,1},y_{20,3}+y_{20,2}],$\\
$D_{46}=[y_{2,4}+y_{2,2},2y_{6,2},y_{16,3}-y_{16,1}],$ &
$D_{47}=[y_{2,3}-y_{2,1},y_{16,4}+y_{16,1},y_{24,4}+y_{24,1}],$\\
$D_{48}=[y_{4,4}-y_{4,2},y_{10,4}+y_{10,2},y_{24,4}+y_{24,2}],$ &
$D_{49}=[y_{14,4}+y_{14,2},2y_{20,2},y_{24,3}-y_{24,2}],$\\
$D_{50}=[y_{2,3}-y_{2,2},y_{4,2}+y_{4,1},y_{16,4}-y_{16,1}],$ &
$D_{51}=[y_{2,2}+y_{2,1},2y_{10,2},y_{22,4}-y_{22,2}],$\\
$D_{52}=[y_{1,3}-y_{1,2},y_{1,3}+y_{1,2},y_{8,4}+y_{8,2}],$ &
$D_{53}=[y_{4,2}-y_{4,1},y_{6,2}-y_{6,1},y_{24,3}-y_{24,1}],$\\
$D_{54}=[y_{1,3}-y_{1,1},y_{1,3}+y_{1,1},y_{6,4}+y_{6,2}],$ &
$D_{55}=[y_{6,3}-y_{6,1},y_{8,3}-y_{8,1},y_{20,2}-y_{20,1}],$\\
$D_{56}=[y_{2,4}-y_{2,2},y_{8,4}-y_{8,2},y_{14,3}-y_{14,2}],$ &
$D_{57}=[y_{2,3}+y_{2,1},y_{18,4}-y_{18,1},y_{20,2}+y_{20,1}],$\\
$D_{58}=[y_{6,2}+y_{6,1},y_{6,3}-y_{6,2},y_{14,4}-y_{14,2}],$ &
$D_{59}=[y_{4,2}+y_{4,1},y_{18,4}+y_{18,1},y_{20,4}+y_{20,2}],$\\
$D_{60}=[y_{6,3}+y_{6,1},y_{14,3}+y_{14,2},y_{22,4}-y_{22,3}],$ &
$D_{61}=[y_{18,4}-y_{18,2},y_{22,4}+y_{22,2},y_{24,3}+y_{24,1}],$\\
$D_{62}=[y_{4,2}-y_{4,1},y_{8,3}+y_{8,1},y_{24,4}-y_{24,2}].$\\
\end{tabular}}}\caption{$D_l$, $0\leq l\leq 62$}\label{tab3}
\end{table}

\begin{theorem}\label{thm:125p}
There exists a $(\mathbb{Z}_{125}\times \mathbb{F}_{q}^+,\mathbb{Z}_{125}\times \{0\},6,1)$-FDF for any prime $q\equiv 7 \pmod{12}$ and $q>43$.
\end{theorem}

\Proof Since $q\equiv 7 \pmod{12}$, $-1\in C^{6,q}_3$. By Lemma \ref{lem:125p} (1)-(4), if one can choose an appropriate mapping $g$ acting on symbolic expressions satisfying that
\begin{itemize}
\item $\{g(y_{11}),g(y_{12}),g(y_{13})\}=\{0,1,2\}$,
\item $g(y_{2i,1})\neq g(y_{2i,2})$ for $1\leq i\leq 12$,
\item $\{g(d):\ d\in D_l\}$ is $\{0,1,2\}$ for each $1\leq l\leq 62$,
\end{itemize}
and can choose appropriate elements of $C_j$, $1\leq j\leq 25$, such that these elements are consistent with the mapping $g$, i.e., $g$ can be seen as a function from $\mathbb{F}_q^*$ to $\mathbb{Z}_3$ satisfying $g(x)=\theta$ if $x\in C^{3,q}_\theta$, then one can apply Lemma \ref{125FDF} to obtain a $(\mathbb{Z}_{125}\times \mathbb{F}_{q}^+,\mathbb{Z}_{125}\times \{0\},6,1)$-FDF. Note that once the above second condition is satisfied, let $\alpha_i=\mathbb{Z}_3\setminus\{f(y_{2i,1}),f(y_{2i,2})\}$, $1\leq i\leq 12$, and we require $y_{2i,3}$ and $y_{2i,4}$ belong to different cosets $C^{6,q}_{\alpha_i}$ and $C^{6,q}_{\alpha_i+3}$.

By Lemma \ref{lem:125p} (6) and (7), the key to pick up an appropriate mapping $g$ is to assign values of $g$ for elements of Types (IV) and (V) (note that elements of Type (V) is contained in a unique $D_l$ for some $1\leq l\leq 62$ and is related with some $C_{2i}$,  $1\leq i\leq 12$). We here give explicit values of $g$ for elements of Types (IV) and (V) in Table \ref{tab2}. Then combining Table \ref{tab3}, one can give values of $g$ for elements of Types (I), (II) and (III). For example, for $D_{33}=[y_{10,2}-y_{10,1},y_{18,3}+y_{18,1},y_{18,4}-y_{18,3}]$, by Table \ref{tab2}, $g(y_{10,2}-y_{10,1})=1$, so it suffices to require $\{g(y_{18,3}+y_{18,1}),g(y_{18,4}-y_{18,3})\}=\{0,2\}$.

\begin{table}[h]\centering
\begin{tabular}{|c|c|c|c|c|c|c|c|c|c|c|c|c|c|}
\hline $i$  & $1$ & $2$ & $3$ & $4$ & $5$ & $6$ & $7$ & $8$ & $9$ & $10$ & $11$ & $12$\\
\hline $g(2y_{2i,1})$  & $1$ & $1$ & $1$ & $2$ & $1$ & $2$ & $0$ & $0$ & $1$ & $2$ & $0$ & $2$\\
\hline $g(2y_{2i,2})$  & $2$ & $2$ & $2$ & $1$ & $2$ & $1$ & $2$ & $1$ & $0$ & $0$ & $2$ & $0$\\
\hline $g(y_{2i,2}-y_{2i,1})$  & $2$ & $1$ & $2$ & $1$ & $1$ & $2$ & $2$ & $2$ & $0$ & $1$ & $1$ & $0$\\
\hline $g(y_{2i,2}+y_{2i,1})$  & $1$ & $1$ & $1$ & $1$ & $1$ & $0$ & $2$ & $1$ & $2$ & $0$ & $0$ & $1$\\
\hline
\end{tabular}\caption{the values of $g$ for elements of Types (IV) and (V) }\label{tab2}
\end{table}

Once $g$ is fixed, one can apply Theorem \ref{thm:cyclot bound} and Lemma \ref{125FDF} to obtain a $(\mathbb{Z}_{125}\times \mathbb{F}_{q}^+,\mathbb{Z}_{125}\times \{0\},6,1)$-FDF for any prime $q\equiv 7 \pmod{12}$ and $q>Q(3,7)=6.43306\times 10^7$. Note that the number $7$ in $Q(3,7)$ is from the fact that the number of cyclotomic conditions on $y_{2i,s}$ is $7$ for $1\leq i\leq 12$ and $s\in\{1,2\}$; the number of cyclotomic conditions on $y_{2i,r}$ is $6$, $1\leq i\leq 12$, $r\in\{3,4\}$; the number of cyclotomic conditions on $y_{1,t}$ is $5$, $t\in\{1,2,3\}$.

For primes $q\equiv 7 \pmod{12}$ and $50023 \leq q\leq Q(3,7)$, by computer search, we can pick up appropriate elements of $C_j$, $1\leq j\leq 25$, such that they are consistent with the mapping $g$ given in Table \ref{tab2}.

On the other hand, elements of $C_j$ do not have to be what they look like in Lemma \ref{125FDF}. It's enough to require each $C_j$, $1\leq j\leq 25$, is a representative system for the cosets of $C_0^{6,q}$ in $\mathbb{F}_q^*$, such that each $\Delta_l$, $l\in \mathbb{Z}_{125}$, is also a representative system for the cosets of $C_0^{6,q}$ in $\mathbb{F}_q^*$. If we allow $C_j$ to vary among the multisets that satisfy the required conditions, by the use of computer, we can also find appropriate elements of $C_j$, $1\leq j\leq 25$, for primes $q\equiv 7 \pmod{12}$ and $43<q<50023$. For example for $p=67$, we can take
\begin{center}
\begin{tabular}{ll}
$C_1=\{ 1, -1, 6, -6, 7, -7 \}$,& $C_2=\{ 1, -1, 2, -2, 4, 20 \}$,\\
$C_4=\{ 1, -1, 2, -2, 4, 11 \}$,&  $C_6=\{ 1, -1, 17, -17, 12, 29 \}$,\\
$C_8=\{ 2, -2, 1, -1, 4, 32 \}$,& $C_{10}=\{ 1, -1, 30, -30, 12, 35 \}$,\\
$C_{12}=\{ 2, -2, 5, -5, 20, 4 \}$,&  $C_{14}=\{ 1, 43, 13, 19, 4, 46 \}$,\\
\end{tabular}
\end{center}

\begin{center}
\begin{tabular}{ll}
$C_{16}=\{ 1, 5, 13, 16, 31, 36 \}$,& $C_{18}=\{ 1, 3, 10, 44, 46, 33 \}$,\\
$C_{20}=\{ 1, 53, 17, 50, 63, 21 \}$,&  $C_{22}=\{ 2, 37, 63, 4, 42, 9 \}$,\\
$C_{24}=\{ 2, 6, 53, 1, 35, 12 \}$.
\end{tabular}
\end{center}
and $C_{2i+1}=-C_{2i}$ for $1\leq i\leq 12$. The interested reader may get a copy of these data from the authors. \qed

\begin{lemma}\label{lem:5p}{\rm (Theorem $4.4$ in \cite{bf})}
There exists a $(\mathbb{Z}_{5}\times \mathbb{F}_{q}^+,\mathbb{Z}_{5}\times \{0\},6,1)$-FDF for any prime $q\equiv 1 \pmod{12}$ and $q>37$.
\end{lemma}

\begin{theorem}\label{125FDF by recur}
There exists a $(\mathbb{F}_{125}^+\times \mathbb{F}_{q}^+,\mathbb{F}_{125}^+\times \{0\},6,1)$-FDF for any prime $q\equiv 1 \pmod{12}$ and $q>37$.
\end{theorem}

\Proof By Lemma \ref{lem:5p}, there exists a $(\mathbb{Z}_{5}\times \mathbb{F}_{q}^+,\mathbb{Z}_{5}\times \{0\},6,1)$-FDF for any prime $q\equiv 1 \pmod{12}$ and $q>37$. Start from this FDF and then apply Construction \ref{recursive} with a homogeneous $(\mathbb{F}_{25}^+,6,1)$-DM to obtain a $(\mathbb{Z}_{5}\times \mathbb{F}_{q}^+\times \mathbb{F}_{25}^+,\mathbb{Z}_{5}\times \{0\}\times \mathbb{F}_{25}^+,6,1)$-FDF. Note that $\mathbb{Z}_{5}\times \mathbb{F}_{25}^+$ is isomorphic to $\mathbb{F}_{125}^+$. \qed

\begin{theorem}\label{6,1RBIBD2}
\begin{itemize}
\item[$(1)$] There exists a $(125q+1,6,1)$-RBIBD for any prime $q\equiv 7 \pmod{12}$ and $q>43$.
\item[$(2)$] There exists a 1-rotational $(125q+1,6,1)$-RBIBD over $\mathbb{F}_{125}^+\times \mathbb{F}_{q}^+$ for any prime $q\equiv 1 \pmod{12}$ and $q>37$.
\end{itemize}
\end{theorem}

\Proof (1) For any prime $q\equiv 7 \pmod{12}$ and $q>43$, by Theorem \ref{thm:125p}, there exists a $(\mathbb{Z}_{125}\times \mathbb{F}_{q}^+,\mathbb{Z}_{125}\times \{0\},6,1)$-FDF. Then apply Proposition \ref{FDFtoRBIBD} with a $(126,6,1)$-RBIBD (a unital design (cf. \cite{b})) to get a $(125q+1,6,1)$-RBIBD.

(2) By Theorem \ref{125FDF by recur}, there exists a $(\mathbb{F}_{125}^+\times \mathbb{F}_{q}^+,\mathbb{F}_{125}^+\times \{0\},6,1)$-FDF for any prime $q\equiv 1 \pmod{12}$ and $q>37$. Then apply Proposition \ref{FDFtoRBIBD 1-rotational} with a 1-rotational $(126,6,1)$-RBIBD over $\mathbb{F}_{125}^+$, which exists by Example 16.92 in \cite{ab}, to get a 1-rotational $(125q+1,6,1)$-RBIBD over $\mathbb{F}_{125}^+\times \mathbb{F}_{q}^+$. \qed

\section{Applications}

In this section we establish asymptotic existences for optimal constant composition codes and strictly optimal frequency hopping sequences by the use of frame difference families obtained in this paper.

\subsection{Partitioned difference families and constant composition codes}

Let $(G,+)$ be an abelian group of order $g$ with a subgroup $N$ of order $n$. A $(G,N,K,\lambda)$ {\em partitioned relative difference family} $(PRDF)$ is a family $\mathfrak{B}=[B_1,B_2,\dots,B_r]$ of $G$ such that the elements of $\mathfrak{B}$ form a partition of $G\setminus N$, and the list
$$\Delta \mathfrak{B}:=\bigcup_{i=1}^r[x-y:x,y\in B_i, x\not=y]=\underline{\lambda}(G\setminus N),$$
where $K$ is the multiset $\{|B_i|:1\leq i\leq r\}$. When $N=\{0\}$, a $(G,\{0\},K,\lambda)$-PRDF is called a {\em partitioned difference family} and simply written as a $(G,K,\lambda)$-PDF. The members of $\mathfrak{B}$ are called {\em base blocks}.

We often use an exponential notation to describe the multiset $K$: a $(G,N,[k_1^{u_1} k_2^{u_2}$ $\cdots k_l^{u_l}],\lambda)$-PRDF is a PRDF in which there are $u_j$ base blocks of size $k_j$, $1\leq j\leq l$.

\begin{proposition}\label{FDFtoPDF}
If there exists an elementary $(G,N,k,1)$-FDF with $|N|=k-1$, then there exists a $(G,[(k-1)^1 k^{s}],k-1)$-PDF, where $s=(|G|-k+1)/k$.
\end{proposition}

\Proof Let $\mathfrak{F}$ be an elementary $(G,N,k,1)$-FDF with $|N|=k-1$. Then $\mathfrak{F}$ satisfies $\bigcup_{F\in \mathfrak{F}, h\in N}(F+h)=G\setminus N$. Set $$\mathfrak{B}=\{F+h:F\in \mathfrak{F}, h\in N\}\cup \{N\}.$$
Then $\mathfrak{B}$ forms a $(G,[(k-1)^1 k^{s}],k-1)$-PDF, where $s=(|G|-k+1)/k$. \qed

Combining the results of Lemma \ref{FDF1}, Theorems \ref{thm:DF-3}, \ref{thm:Twin DS-2}, \ref{thm:Singer DS-2}, \ref{thm:DF-4} and Proposition \ref{FDFtoPDF}, we have

\begin{theorem}\label{cor:PDF}
\begin{itemize}
\item[(1)] There exists a $(\mathbb{Z}_{7}\times \mathbb{F}_{89}^+,[7^1 8^{77}],7)$-PDF.

\item[(2)] There exists an $(\mathbb{F}_{p}^+ \times \mathbb{F}_{q}^+,[p^1 (p+1)^{s}],p)$-PDF
\begin{itemize}
\item for any prime power $p\equiv 3\pmod{4}$ and any prime power $q\equiv 1 \pmod{p+1}$ with $q>Q((p+1)/2,p)$;
\item for any prime power $p\equiv 1\pmod{4}$ and any prime power $q\equiv 1 \pmod{2p+2}$ with $q>Q(p+1,p)$,
\end{itemize}
where $s=p(q-1)/(p+1)$.

\item[(3)] Let $p$ and $p+2$ be twin prime powers satisfying $p>2$. There exists an $(\mathbb{F}_{p}^+ \times \mathbb{F}_{p+2}^+\times\mathbb{F}_{q}^+,[(p(p+2))^1 (p(p+2)+1)^{s}],p(p+2))$-PDF for any prime power $q\equiv 1 \pmod{p(p+2)+1}$ and $q>Q((p(p+2)+1)/2,p(p+2))$, where $s=p(p+2)(q-1)/(p(p+2)+1)$.

\item[(4)] There exists a $(\mathbb{Z}_{2^m-1} \times \mathbb{F}_{q}^+,[(2^m-1)^1 (2^m)^{s}],2^m-1)$-PDF for any integer $m\geq3$ and any prime power $q\equiv 1 \pmod{2^{m}}$ with $q>Q(2^{m-1},2^{m}-1)$, where $s=(2^m-1)(q-1)/2^m$.
\end{itemize}
\end{theorem}

Partitioned difference families were explicitly introduced in \cite{dy} to construct optimal constant composition codes, which can be used in the MFSK modulation of power line communications (cf. \cite{pvyh}). For more information on partitioned difference families and its relationship with other topics such as zero-difference balanced functions (cf. \cite{zty}), the interested reader may refer to \cite{lwg}. We remark that recently M. Buratti \cite{b17} established a construction for partitioned difference families by the use of Hadamard strong difference families.

Let $Q=\{0,1,\ldots,q-1\}$ be an alphabet with $q$ symbols. An $(n,M,d,[\omega_0$, $\omega_1,\ldots,$ $\omega_{q-1}])_q$ constant composition code $($CCC$)$
is a subset ${\cal C}\subseteq Q^n$ with size $M$ and minimum Hamming distance $d$, such that the symbol $i$, $i\in Q$, appears exactly $\omega_i$ times in each codeword of $\cal C$.

Since there is no essential difference among the symbols of $Q$, one often regards $[\omega_0,\omega_1,\ldots,\omega_{q-1}]$ as a multiset and denotes it by $[k_1^{u_1}k_2^{u_2}\cdots k_l^{u_l}]$, where $k_j$ appears $u_j$ times in $[\omega_0,\omega_1,\ldots,\omega_{q-1}]$, $1\leq j\leq l$. Let $A_q(n,d,[\omega_0,\omega_1,\ldots,\omega_{q-1}])$ be the maximal size of an $(n,M, d,[\omega_0,\omega_1,\ldots,\omega_{q-1}])_q$-CCC.

\begin{proposition}\label{prop:CCC}{\rm \cite{lvc}}
If $nd-n^2+(\omega_0^2+\omega_1^2+\cdots+\omega_{q-1}^2)>0$, then
\begin{eqnarray}\label{ccc-bound}
A_q(n,d,[\omega_0,\omega_1,\ldots,\omega_{q-1}])\leq \frac{nd}{nd-n^2+(\omega_0^2+\omega_1^2+\cdots+\omega_{q-1}^2)}.
\end{eqnarray}
\end{proposition}

A constant composition code attaining the bound (\ref{ccc-bound}) is called {\em optimal}. The following proposition indicates that optimal CCCs can be derived from PDFs.

\begin{proposition}\label{prop:CCC-PDF-relation}{\rm (Construction $6$ in \cite{dy})}
If a $(G,[k_1^{u_1}k_2^{u_2}\cdots k_l^{u_l}],\lambda)$-PDF exists, then there is an optimal $(|G|,|G|,|G|-\lambda,[k_1^{u_1}k_2^{u_2}\cdots k_l^{u_l}])_q$-CCC meeting the bound (\ref{ccc-bound}), where $q=\sum_{j=1}^l u_j$.
\end{proposition}

Applying Proposition \ref{prop:CCC-PDF-relation} with PDFs from Theorem \ref{cor:PDF}, we have

\begin{theorem}\label{thm:ccc}
\begin{itemize}
\item[(1)] There exists an optimal $(623,623,616,[7^1 8^{77}])_{78}$-CCC.
\item[(2)] There exists an optimal $(pq,pq,p(q-1),[p^1 (p+1)^{s}])_{s+1}$-CCC
\begin{itemize}
\item for any prime power $p\equiv 3\pmod{4}$ and any prime power $q\equiv 1 \pmod{p+1}$ with $q>Q((p+1)/2,p)$;
\item for any prime power $p\equiv 1\pmod{4}$ and any prime power $q\equiv 1 \pmod{2p+2}$ with $q>Q(p+1,p)$,
\end{itemize}
where $s=p(q-1)/(p+1)$.

\item[(3)] Let $p$ and $p+2$ be twin prime powers satisfying $p>2$. There exists an optimal $(p(p+2)q,p(p+2)q,p(p+2)(q-1),[(p(p+2))^1 (p(p+2)+1)^{s}])_{s+1}$-PDF for any prime power $q\equiv 1 \pmod{p(p+2)+1}$ and $q>Q((p(p+2)+1)/2,p(p+2))$, where $s=p(p+2)(q-1)/(p(p+2)+1)$.

\item[(4)] There exists an optimal $((2^m-1)q,(2^m-1)q,(2^m-1)(q-1),[(2^m-1)^1 (2^m)^{s}])_{s+1}$-PDF for any integer $m\geq3$ and any prime power $q\equiv 1 \pmod{2^{m}}$ with $q>Q(2^{m-1},2^{m}-1)$, where $s=(2^m-1)(q-1)/2^m$.
\end{itemize}
\end{theorem}

\subsection{Frequency hopping sequences}

Frequency hopping multiple-access has been widely used in the modern communication systems such as
ultrawideband, military communications and so on (cf. \cite{pd,yb}).

Let $F=\{f_0,f_1,\ldots,f_{l-1}\}$ be a set $($called an {\em alphabet}$)$ of $l\geq 2$ available frequencies. A sequence $X=\{x(t)\}_{t=0}^{n-1}$ is called a {\em frequency hopping sequence} $($FHS$)$ of length $n$ over $F$ if $x(t)\in F$ for any $0\leq t\leq n-1$.

For any FHS $X=\{x(t)\}_{t=0}^{n-1}$, {\em the partial Hamming autocorrelation} function of X for a correlation window length $L$ starting at $j$ is defined by
\begin{eqnarray}\label{autocorrelation}
H_{X,X}(\tau;j|L)=\sum_{t=j}^{j+L-1}h[x(t),x(t+\tau)],0\leq \tau\leq n-1,
\end{eqnarray}
where $1\leq L\leq n$, $0\leq j\leq n-1$, $h[a,b]=1$ if $a=b$ and 0 otherwise, and the addition is performed modulo $n$. If $L=n$, the partial Hamming correlation function defined in (\ref{autocorrelation}) becomes the {\em conventional periodic Hamming correlation} (cf. \cite{lg}).

For any FHS $X=\{x(t)\}_{t=0}^{n-1}$ and any given $1\leq L\leq n$, define
$$H(X;L)=\max_{0\leq j<n} \max_{1\leq \tau <n}\{H_{X,X}(\tau;j|L)\}.$$

\begin{proposition}\label{prop:FHS}{\rm \cite{czyt}}
Let $X$ be an FHS of length $n$ over an alphabet of size $l$. Then, for each window length $L$ with $1\leq L\leq n$,
\begin{eqnarray}\label{FHS-bound}
H(X;L)\geq \left\lceil\frac{L}{n}\left\lceil\frac{(n-\epsilon)(n+\epsilon-l)}{l(n-1)}\right\rceil\right\rceil,
\end{eqnarray}
where $\epsilon$ is the least nonnegative residue of $n$ modulo $l$.
\end{proposition}

Let $X$ be an FHS of length $n$ over an alphabet $F$. It is said to be {\em strictly optimal} if the bound $(\ref{FHS-bound})$ in Proposition \ref{prop:FHS} is met for any $1\leq L\leq n$.

\begin{proposition}{\rm (Theorem $3.7$ in \cite{bj})}
Let $k$ and $v$ be positive integers satisfying $k+1|v-1$. Then there exists a strictly optimal FHS of  length $kv$ over an alphabet of size $(kv+1)/(k+1)$ if and only if there exists an elementary $(kv,k,k+1,1)$-FDF over $\mathbb{Z}_{kv}$.
\end{proposition}

Combining the results of Lemma \ref{FDF1}, Theorems \ref{thm:DF-3}, \ref{thm:Twin DS-2}, \ref{thm:Singer DS-2}, \ref{thm:DF-4} and Proposition \ref{prop:FHS}, we have the following theorem. Note that to apply Proposition \ref{prop:FHS}, the needed FDFs must be defined on a cyclic group.

\begin{theorem}\label{thm:FHS}
\begin{itemize}
\item[(1)] There exists a strictly optimal FHS of length $623$ over an alphabet of size $78$.
\item[(2)] There exists a strictly optimal FHS of length $pq$ over an alphabet of size $(pq+1)/(p+1)$
\begin{itemize}
\item for any prime $p\equiv 3\pmod{4}$ and any prime $q\equiv 1 \pmod{p+1}$ with $q>Q((p+1)/2,p)$;
\item for any prime $p\equiv 1\pmod{4}$ and any prime $q\equiv 1 \pmod{2p+2}$ with $q>Q(p+1,p)$.
\end{itemize}

\item[(3)] Let $p$ and $p+2$ be twin primes. There exists a strictly optimal FHS of length $p(p+2)q$ over an alphabet of size $(p(p+2)q+1)/(p(p+2)+1)$ for any prime $q\equiv 1 \pmod{p(p+2)+1}$ and $q>Q((p(p+2)+1)/2,p(p+2))$.

\item[(4)] There exists a strictly optimal FHS of length $(2^m-1)q$ over an alphabet of size $((2^m-1)q+1)/2^m$ for any integer $m\geq3$ and any prime $q\equiv 1 \pmod{2^{m}}$ with $q>Q(2^{m-1},2^{m}-1)$.
\end{itemize}
\end{theorem}

\section{Concluding remarks}

By a careful application of cyclotomic conditions attached to strong difference families, this paper establishes (asymptotic) existences of several classes of frame difference families, which are used to derive new resolvable balanced incomplete block designs, new optimal constant composition codes and new strictly optimal frequency hopping sequences.

We believe that starting from those CCCs or FHSs obtained in Section 6, and applying appropriate known recursive constructions in the literature, one can obtain more new existence results on them. For example, by Construction 3 in \cite{lwg} and via similar technique in the proof of Theorems 18 and 19 in \cite{lwg}, one can obtain more new PDFs, which can yield new CCCs; apply Theorem 6.8 in \cite{bj} to obtain more new FHSs, and so on.

Frame difference families can be seen as special resolvable difference families (cf. \cite{b97,byw}). By using a $(44,4,5,2)$ resolvable difference family, which is not a frame difference family, M. Buratti, J. Yan and C. Wang \cite{byw} presented the first example of a $(45,5,2)$-RBIBD. Thus an interesting future direction is to establish constructions, especially direct constructions, for resolvable difference families.

Most of our results in this paper rely heavily on existences of elements satisfying certain cyclotomic conditions in a finite field, and Theorem \ref{thm:cyclot bound} just supplies us a way to ensure existences of such elements. We remark that recently X. Lu improved the lower bound on $q$ in Theorem \ref{thm:cyclot bound} in some circumstances (see Theorem 3 in \cite{l}). He introduced an existence bound $L(d,t)$ for elements $x$ that satisfies the following cyclotomic conditions:
\begin{itemize}
\item[i)] $x\in \bigcup_{i\in U(\mathbb{Z}_d)}\  C_i^{d, q}$, where $U(\mathbb{Z}_d)$ is the set of all units in $\mathbb{Z}_d$;
\item[ii)] $x^{d-c_j}(a_jx+b_j)\in C_{0}^{d, q}$ for $1\leq j\leq t-1$.
\end{itemize}
One possible development of this paper could be to find suitable constructions that use his result.

\subsection*{Acknowledgements}
The authors would like to thank Professor Marco Buratti of Universit\`a di Perugia for his many valuable comments.

\end{document}